# On the existence of regular vectors

Christoph Zellner[*]

**Abstract.** Let $G$ be a locally convex Lie group and $\pi : G \to \mathrm{U}(\mathcal{H})$ be a continuous unitary representation. $\pi$ is called smooth if the space of $\pi$-smooth vectors $\mathcal{H}^\infty \subset \mathcal{H}$ is dense. In this article we show that under certain conditions, concerning in particular the structure of the Lie algebra $\mathfrak{g}$ of $G$, a continuous unitary representation of $G$ is automatically smooth. As an application, this yields a dense space of smooth vectors for continuous positive energy representations of oscillator groups, double extensions of loop groups and the Virasoro group. Moreover we show the existence of a dense space of analytic vectors for the class of semibounded representations of Banach–Lie groups. Here $\pi$ is called semibounded, if $\pi$ is smooth and there exists a non-empty open subset $U \subset \mathfrak{g}$ such that the operators $i\mathrm{d}\pi(x)$ from the derived representation are uniformly bounded from above for $x \in U$.

**2010 Mathematics Subject Classification.** Primary 22E65; Secondary 22E45

**Keywords.** infinite dimensional Lie group, smooth vector, analytic vector, semibounded representation, continuous representation

## 1. Introduction

Let $G$ be a Lie group modeled on a locally convex space and $\mathfrak{g}$ be its Lie algebra, cf. [11] for the basic concepts of infinite-dimensional Lie theory. Assume that $G$ has an *exponential map*, i.e., a smooth map $\exp : \mathfrak{g} \to G$ such that, for every $x \in \mathfrak{g}$, the curve $\gamma_x(t) := \exp(tx)$ is a one-parameter group with $\gamma'_x(0) = x$. Let $\pi : G \to \mathrm{U}(\mathcal{H})$ be a *continuous unitary representation* of $G$, i.e., a homomorphism $\pi$ of $G$ into the unitary group of a complex Hilbert space $\mathcal{H}$ such that $\pi$ is continuous when $\mathrm{U}(\mathcal{H})$ is equipped with the strong operator topology.

A vector $v \in \mathcal{H}$ is called a *smooth vector* if the orbit map $\pi^v : G \to \mathcal{H}, \pi^v(g) := \pi(g)v$ is smooth and we denote the space of all smooth vectors by $\mathcal{H}^\infty = \mathcal{H}^\infty(\pi)$. Then $\mathcal{H}^\infty$ is a $\pi$-invariant subspace of $\mathcal{H}$. We call $\pi$ *smooth* if $\mathcal{H}^\infty \subset \mathcal{H}$ is dense. In contrast to finite-dimensional Lie groups not every continuous unitary representation $\pi$ of an infinite-dimensional Lie group $G$ is smooth (see [14]). In this article we obtain a result which states that under certain conditions a continuous unitary representation is automatically smooth. In particular, this applies to positive energy representations of some prominent groups arising in physics, which will be considered in detail. If $G$ is a Banach–Lie group a vector $v \in \mathcal{H}$ is called an *analytic vector* if $\pi^v : G \to \mathcal{H}$ is analytic. We will show the existence of a dense space of analytic vectors for the class of so-called semibounded representations of Banach–Lie groups. Semiboundedness is a stable version of the 'positive energy' condition ([1, 2]) and it will be defined below.

---
[*]The author acknowledges the support of DFG-grant NE 413/7-2 in the framework of SPP "Representation Theory". The author is grateful to Karl-Hermann Neeb for helpful discussions and for reading a preliminary version of the manuscript.

For a continuous representation $\pi : G \to \mathrm{U}(\mathcal{H})$ we obtain for every $x \in \mathfrak{g}$ a continuous unitary one-parameter group $\pi_x(t) := \pi(\exp(tx))$. By Stone's Theorem $\pi_x$ has a skew-adjoint generator $\overline{\mathrm{d}\pi}(x)$ with dense domain

$$\mathcal{D}(\overline{\mathrm{d}\pi}(x)) \subset \mathcal{H}, \quad \pi(\exp(tx)) = e^{t\overline{\mathrm{d}\pi}(x)}.$$

The operators $i\overline{\mathrm{d}\pi}(x)$ are self-adjoint and we define the support functional

$$s_\pi : \mathfrak{g} \to \mathbb{R} \cup \{\infty\}, \quad s_\pi(x) := \sup(\mathrm{Spec}(i\overline{\mathrm{d}\pi}(x))),$$

which is easily seen to be invariant under the adjoint action $\mathrm{Ad}$ of $G$ on $\mathfrak{g}$. Moreover

$$W_\pi := \{x_0 \in \mathfrak{g} : s_\pi \text{ is bounded in a neighborhood of } x_0\}$$

is a (possibly empty) invariant cone in $\mathfrak{g}$. In general, the unbounded operators $\overline{\mathrm{d}\pi}(x), x \in \mathfrak{g}$, may not have a common dense domain. However, if $\pi$ is smooth, then every $\overline{\mathrm{d}\pi}(x), x \in \mathfrak{g}$, is essentially skew-adjoint on $\mathcal{H}^\infty$ and we obtain the *derived representation*

$$\mathrm{d}\pi : \mathfrak{g} \to \mathrm{End}(\mathcal{H}^\infty), \quad \mathrm{d}\pi(x) := \overline{\mathrm{d}\pi}(x)|_{\mathcal{H}^\infty},$$

which is a representation of $\mathfrak{g}$ on $\mathcal{H}^\infty$, see [13]. Moreover for a smooth representation $\pi$ the support functional satisfies

$$s_\pi(x) = \sup_{v \in \mathcal{H}^\infty, \|v\|=1} \langle i\mathrm{d}\pi(x)v, v \rangle, \quad x \in \mathfrak{g} \tag{1}$$

so that the cone $W_\pi$ is convex and $s_\pi$ is lower semicontinuous for $\pi$ smooth ([13]). The representation $\pi$ is called *semibounded* if $\pi$ is smooth and $W_\pi \neq \emptyset$.

It is an open problem if every continuous representation $\pi$ with $W_\pi \neq \emptyset$ is already semibounded. In Section 2 we show that this is true under certain conditions involving conditions on the structure of the Lie algebra $\mathfrak{g}$ of $G$. This result will be applied to continuous positive energy representations of oscillator groups, double extensions of loop groups and the Virasoro group, yielding a dense space of smooth vectors for those. Section 3 is devoted to the existence of analytic vectors. We show that every semibounded representation $\pi : G \to \mathrm{U}(\mathcal{H})$ of a Banach–Lie group $G$ extends locally to a holomorphic map which in turn yields a dense subspace of analytic vectors. This generalizes a result in [10] where it is assumed that an Olshanski semigroup $G \exp(iW_\pi) \subset G_\mathbb{C}$ in a complexification $G_\mathbb{C}$ of $G$ exists, cf. [10, Def. 5.3 & Thm. 5.7]. Up to this result for semibounded representations of Banach–Lie groups, analytic vectors were known to exist only for very special classes of groups such as certain direct limits ([23]) and the canonical commutation relations in Quantum Field Theory ([21]).

## 2. Smooth vectors

In this section $G$ denotes a (locally convex) Lie group with an exponential map $\exp : \mathfrak{g} \to G$. For a continuous representation $\pi : G \to \mathrm{U}(\mathcal{H})$ and a subset $\mathfrak{h} \subset \mathfrak{g}$ we define inductively

$$\mathcal{D}(\mathfrak{h}) := \mathcal{D}^1(\mathfrak{h}) := \bigcap_{x \in \mathfrak{h}} \mathcal{D}(\overline{\mathrm{d}\pi}(x)),$$

$$\mathcal{D}^n(\mathfrak{h}) := \{v \in \mathcal{D}(\mathfrak{h}) : (\forall x \in \mathfrak{h})\ \overline{\mathrm{d}\pi}(x)v \in \mathcal{D}^{n-1}(\mathfrak{h})\},$$

$$\mathcal{D}^\infty(\mathfrak{h}) := \bigcap_{n \in \mathbb{N}} \mathcal{D}^n(\mathfrak{h}),$$

$\mathcal{D}^n(\overline{\mathrm{d}\pi}(x)) := \mathcal{D}^n(\{x\}), x \in \mathfrak{g}, n \in \mathbb{N} \cup \{\infty\}$ and set

$$\omega_v^k : \mathfrak{g}^k \to \mathcal{H}, \quad \omega_v^k(x_1, \ldots x_k) := \overline{\mathrm{d}\pi}(x_1) \cdots \overline{\mathrm{d}\pi}(x_k) v$$

for $v \in \mathcal{D}^n(\mathfrak{g}), k \leq n$. Furthermore we put

$$\mathcal{D}_c^n(\mathfrak{g}) := \{v \in \mathcal{D}^n(\mathfrak{g}) : \omega_v^k \text{ is continuous and } k\text{-linear for all } k \leq n\}$$

for $n \in \mathbb{N} \cup \{\infty\}$.

### 2.1. General results.

**Definition 2.1.** Assume that $G$ is a Fréchet–Lie group. Denote by $L^\infty([0,1], \mathfrak{g})$ the space of measurable maps $[0,1] \to \mathfrak{g}$ with bounded, separable image modulo those maps vanishing almost everywhere. Equip $L^\infty([0,1], \mathfrak{g})$ with the topology generated by the seminorms $\|\xi\|_{L^\infty, p} := \mathrm{ess\,sup}(p \circ \xi)$ where $p$ runs through the continuous seminorms on $\mathfrak{g}$. Let $R([0,1], \mathfrak{g})$ be the closure of the space of Riemannian step functions $[0,1] \to \mathfrak{g}$ in $L^\infty([0,1], \mathfrak{g})$. Then the Lie group $G$ is called *R-regular* if, for every $\xi \in R([0,1], \mathfrak{g})$, the initial value problem [1]

$$\gamma'(t) = \gamma(t).\xi(t), \quad \gamma(0) = \mathbf{1}$$

has a solution $\gamma_\xi$ on $[0,1]$ (which then is uniquely determined by the product rule) and the corresponding map

$$\mathrm{Evol} : R([0,1], \mathfrak{g}) \to C([0,1], G), \quad \xi \mapsto \gamma_\xi$$

is smooth, cf. [6] for more details.

**Remark 2.2.** Let $\pi : G \to \mathrm{U}(\mathcal{H})$ be a continuous unitary representation. For every $v \in \mathcal{H}$ the map

$$C([0,1], G) \to C([0,1], \mathcal{H}), \quad \gamma \mapsto \pi^v \circ \gamma$$

---
[1] Here $\gamma : [0,1] \to G$ is an absolutely continuous map and $\gamma'$ denotes its (almost everywhere) defined derivative. Moreover $\gamma(t).\xi(t)$ denotes the tangent vector at $\gamma(t)$ obtained from $\xi(t) \in \mathfrak{g}$ by left translation (via the derivative of left multiplication by $\gamma(t)$).

is continuous. Indeed, for $\varepsilon > 0$ there exists an open **1**-neighborhood $U \subset G$ such that $\|\pi(g)v - v\| < \varepsilon$ for all $g \in U$. For an element $\gamma_0$ in the Lie group $C([0,1], G)$ we then have

$$\sup_t \|\pi(\gamma(t))v - \pi(\gamma_0(t))v\| = \sup_t \|\pi(\gamma_0(t)^{-1}\gamma(t))v - v\| < \varepsilon \quad \forall \gamma \in \gamma_0 \cdot C([0,1], U)$$

and $\gamma_0 \cdot C([0,1], U)$ is an open neighborhood of $\gamma_0$ in $C([0,1], G)$. We similarly obtain that, for every $v \in \mathcal{D}_c^1(\mathfrak{g})$, the linear map

$$R([0,1], \mathfrak{g}) \to L^\infty([0,1], \mathcal{H}), \quad \xi \mapsto (t \mapsto \mathtt{d}\pi(\xi(t))v)$$

is continuous.

The following proposition is a version of [14, Lem. 3.4] for $R$-regular Lie groups.

**Proposition 2.3.** *Let $G$ be a $R$-regular Fréchet–Lie group and $\pi : G \to \mathrm{U}(\mathcal{H})$ a continuous unitary representation. Then $\mathcal{H}^\infty = \mathcal{D}_c^\infty(\mathfrak{g})$.*

*Proof.* The inclusion $\mathcal{H}^\infty \subset \mathcal{D}_c^\infty(\mathfrak{g})$ is clear. Now fix $v \in \mathcal{D}_c^\infty(\mathfrak{g})$ and consider the maps

$$F_1 : R([0,1], \mathfrak{g}) \to C([0,1], \mathcal{H}), \quad F_1(\xi)(t) = \pi(\mathrm{Evol}(\xi)(t))v$$

$$F_2 : R([0,1], \mathfrak{g}) \to C([0,1], \mathcal{H}), \quad F_2(\xi)(t) = v + \int_0^t \pi(\mathrm{Evol}(\xi)(s))\overline{\mathtt{d}\pi}(\xi(s))v \, ds$$

Since $G$ is $R$-regular, $F_1$ is continuous by Remark 2.2. Let $\xi_0, \xi \in R([0,1], \mathfrak{g}), t \in [0,1], w \in \mathcal{H}$ and set $\gamma := \mathrm{Evol}(\xi), \gamma_0 := \mathrm{Evol}(\xi_0)$. Then

$$|\langle (F_2(\xi) - F_2(\xi_0))(t), w \rangle| \leq \int_0^t |\langle \pi(\gamma(s))\overline{\mathtt{d}\pi}(\xi(s))v - \pi(\gamma_0(s))\overline{\mathtt{d}\pi}(\xi_0(s))v, w \rangle| ds$$

$$\leq \int_0^t \|\overline{\mathtt{d}\pi}(\xi(s) - \xi_0(s))v\| \cdot \|w\| + \|\overline{\mathtt{d}\pi}(\xi_0(s))v\| \cdot \|\pi(\gamma(s)^{-1})w - \pi(\gamma_0(s)^{-1})w\| ds$$

$$\leq \sup_{s \in [0,1]} \|\overline{\mathtt{d}\pi}(\xi(s) - \xi_0(s))v\| \cdot \|w\| + c \cdot \sup_{s \in [0,1]} \|\pi(\gamma(s)^{-1})w - \pi(\gamma_0(s)^{-1})w\|,$$

where $c := \sup_{s \in [0,1]} \|\overline{\mathtt{d}\pi}(\xi_0(s))v\| < \infty$. Since Evol is continuous we conclude with Remark 2.2 that $F_2$ is continuous as a map $F_2 : R([0,1], \mathfrak{g}) \to C([0,1], \mathcal{H}_w)$, where $\mathcal{H}_w$ denotes $\mathcal{H}$ equipped with the weak topology. If $\xi$ is a constant map $\xi = x, x \in \mathfrak{g}$, then $F_1(\xi)(t) = \pi(\exp(tx))v$ and

$$F_2(\xi)(t) = v + \int_0^t \pi(\exp(sx))\overline{\mathtt{d}\pi}(x)v \, ds = \pi(\exp(tx))v = F_1(\xi)(t).$$

If $\xi$ is a Riemannian step function $\xi|_{[t_j, t_{j+1}[} = x_j \in \mathfrak{g}, 0 = t_0 < \cdots < t_k = 1$, then

$$\mathrm{Evol}(\xi)|_{[t_j, t_{j+1}]}(t) = \prod_{\ell=0}^{j-1} \exp((t_{\ell+1} - t_\ell)x_\ell) \exp((t - t_j)x_j)$$

and one similarly shows that $F_1$ and $F_2$ coincide on Riemannian step functions. As $F_1$ and $F_2$ are both continuous maps $R([0,1],\mathfrak{g}) \to C([0,1],\mathcal{H}_w)$ we conclude that $F_1 = F_2$. Now let $\gamma \in C^\infty([0,1],G), \gamma(0) = \mathbf{1}$, and denote by $\xi := \delta^l(\gamma) \in C^\infty([0,1],\mathfrak{g})$ the left logarithmic derivative of $\gamma$ so that $\mathrm{Evol}(\xi) = \gamma$, cf. [11, Sect. II.4]. Since $F_1 = F_2$ and $C^\infty([0,1],\mathfrak{g}) \subset R([0,1],\mathfrak{g})$ we obtain that $\pi^v \circ \gamma = F_1(\xi) = F_2(\xi)$ is $C^1$ with derivative

$$\tfrac{d}{dt}(\pi^v \circ \gamma)(t) = \tfrac{d}{dt}(F_2(\xi))(t) = \pi(\gamma(t))\overline{\mathrm{d}\pi}(\xi(t))v.$$

Hence [14, Lem. 3.3] yields that, for every $v \in \mathcal{D}_c^\infty(\mathfrak{g})$, the orbit map $\pi^v : G \to \mathcal{H}$ is $C^1$ with derivative

$$d\pi^v(g)(g.x) = \pi(g)\overline{\mathrm{d}\pi}(x)v,$$

where $g \to g.x$ denotes the left invariant vector field on $G$ corresponding to $x \in \mathfrak{g}$. As $\overline{\mathrm{d}\pi}(x)v \in \mathcal{D}_c^\infty(\mathfrak{g})$ we conclude that $d\pi^v$ is $C^1$ and hence $\pi^v$ is $C^2$. Iterating this argument shows that $\pi^v$ is $C^\infty$, i.e., $v \in \mathcal{H}^\infty$. □

**Definition 2.4.** The (locally convex) Lie group $G$ has the *Trotter property* if, for every $x, y \in \mathfrak{g}$,

$$\exp(t(x+y)) = \lim_{n \to \infty} \left( \exp\left(\tfrac{t}{n}x\right) \exp\left(\tfrac{t}{n}y\right) \right)^n$$

holds uniformly on compact subsets on $\mathbb{R}$.

**Remark 2.5.** There are many important classes of Lie groups with the Trotter property, for example locally exponential Lie groups, $\mathrm{Diff}(M)$ for a compact manifold $M$ and direct limits $G = \bigcup_n G_n$ with $G_1 \subset G_2 \subset \ldots$ finite-dimensional Lie groups, see [18, Sect. 3] for a more detailed list. Moreover every $R$-regular Lie group satisfies the Trotter property and many important Lie groups are $R$-regular, cf. [6].

**Lemma 2.6.** *Let $\pi : G \to \mathrm{U}(\mathcal{H})$ be a smooth representation and $\mathfrak{h} \subset \mathfrak{g}$ a dense subspace. Then*

$$W_\pi \cap \mathfrak{h} = \{x_0 \in \mathfrak{h} : s_\pi|_{\mathfrak{h}} \text{ is bounded in a neighborhood of } x_0\}.$$

*Proof.* The inclusion $\subset$ is clear. Now let $x_0 \in \mathfrak{h}$ such that $s_\pi|_{\mathfrak{h}}$ is bounded in a neighborhood of $x_0$, i.e., there exists $U \subset \mathfrak{g}$ open with $x_0 \in U$ and $C > 0$ such that $s_\pi(x) \leq C$ for all $x \in U \cap \mathfrak{h}$. As $s_\pi$ is lower semicontinuous the set $\{a \in \mathfrak{g} : s_\pi(a) \leq C\}$ is closed. Hence $s_\pi(x) \leq C$ for all $x \in \overline{U \cap \mathfrak{h}}$. Since $\mathfrak{h} \subset \mathfrak{g}$ is dense we have $\overline{U \cap \mathfrak{h}} \supset U$. This implies $x_0 \in W_\pi$. □

**Lemma 2.7.** *Let $\pi : G \to \mathrm{U}(\mathcal{H})$ be a continuous representation and $\mathfrak{h} \subset \mathfrak{g}$ a dense subspace. Let $v \in \mathcal{D}(\mathfrak{h})$ such that $\alpha : \mathfrak{h} \to \mathcal{H}, x \mapsto \overline{\mathrm{d}\pi}(x)v$ is continuous and linear. Then $v \in \mathcal{D}(\mathfrak{g})$ and $\overline{\mathrm{d}\pi}(x)v = \widehat{\alpha}(x)$ for all $x \in \mathfrak{g}$, where $\widehat{\alpha} : \mathfrak{g} \to \mathcal{H}$ is the unique continuous linear extension of $\alpha$.*

*Proof.* Since $v \in \mathcal{D}(\mathfrak{h})$ we have

$$\pi(\exp(tx))v - v = \int_0^t \pi(\exp(sx))\widehat{\alpha}(x)v \, ds \qquad (2)$$

for $x \in \mathfrak{h}, t \in \mathbb{R}$. As both sides of (2) are continuous maps in $x \in \mathfrak{g}$ (for each $t \in \mathbb{R}$) we conclude that (2) holds for all $x \in \mathfrak{g}$. This implies $v \in \mathcal{D}(\mathfrak{g})$ and $\overline{d\pi}(x)v = \widehat{\alpha}(x), x \in \mathfrak{g}$. □

**Definition 2.8.** An *integral subgroup* of $G$ is a subgroup $H \subset G$ equipped with a Lie group structure such that $H$ is connected, the inclusion $H \hookrightarrow G$ is a smooth group homomorphism and the induced map of Lie algebras $\mathfrak{h} \to \mathfrak{g}$ is injective, see also [11, Def. IV.4.7].

**Theorem 2.9.** *Let $\pi : G \to \mathrm{U}(\mathcal{H})$ be a continuous representation. Assume there exist subalgebras $\mathfrak{h}_j \subset \mathfrak{g}, j \in J$, with corresponding integral subgroups $H_j \subset G$ such that $\mathfrak{h} := \mathrm{span}_j \mathfrak{h}_j \subset \mathfrak{g}$ is dense and $\pi|_{H_j}$ is smooth for every $j \in J$. Further assume the existence of an element $x_0 \in \bigcap_j \mathfrak{h}_j$ such that $s_\pi|_{\mathfrak{h}}$ is bounded in a neighborhood of $x_0$. Then the following assertions hold:*

(a) *If $G$ has the Trotter property then $\mathcal{D}^\infty(\overline{d\pi}(x_0)) = \mathcal{D}_c^\infty(\mathfrak{g})$.*

(b) *If $G$ is locally exponential then $\pi$ is semibounded.*

*Proof.* By assumption we may choose a continuous seminorm $p$ on $\mathfrak{g}$ and $c > 0$ such that

$$s_\pi(x_0 + y) \leq c \text{ for all } y \in \mathfrak{h}, p(y) \leq 1. \qquad (3)$$

Let $j \in J$. Then $\pi|_{H_j}$ is semibounded since $x_0 \in W_{\pi|_{H_j}}$. Therefore, by [19, Thm. 3.4], $\mathcal{H}^\infty(\pi|_{H_j}) = \mathcal{D}^\infty(\overline{d\pi}(x_0))$ and

$$\|\overline{d\pi}(x)v\| \leq (p(x) + p([x, x_0])) \cdot \|Nv\|, \quad x \in \mathfrak{h}_j, v \in \mathcal{D}^\infty(\overline{d\pi}(x_0)), \qquad (4)$$

where $N := \frac{1}{i}\overline{d\pi}(x_0) + (c+1)\mathbf{1} \geq \mathbf{1}$. Note that $\mathcal{D}^\infty(\overline{d\pi}(x_0)) = \mathcal{D}^\infty(N)$. Since $\pi|_{H_j}$ is smooth, (1) and (3) yield

$$|\langle i\overline{d\pi}(y)v, v\rangle| \leq c\|v\|^2 + \langle \tfrac{1}{i}\overline{d\pi}(x_0)v, v\rangle \leq \langle Nv, v\rangle$$

for $y \in \mathfrak{h}_j, p(y) \leq 1, v \in \mathcal{D}^\infty(\overline{d\pi}(x_0))$. With $y = [x, x_0], x \in \mathfrak{h}_j$ we obtain

$$|\langle [\overline{d\pi}(x), \overline{d\pi}(x_0)]v, v\rangle| \leq p([x, x_0]) \cdot \langle Nv, v\rangle, \quad x \in \mathfrak{h}_j, v \in \mathcal{D}^\infty(\overline{d\pi}(x_0)). \qquad (5)$$

Let $j_1, \ldots, j_n \in J$ and $x_1 \in \mathfrak{h}_{j_1}, \ldots, x_n \in \mathfrak{h}_{j_n}$. Then $\sum_{\ell=1}^n \frac{1}{i}\overline{d\pi}(x_\ell)$ is essentially self-adjoint on $\mathcal{D}^\infty(N)$ by (4), (5) and [22, Thm. X.37]. We now show by induction on $n$ that $\overline{d\pi}(\sum_{\ell=1}^n x_\ell) = \overline{\sum_{\ell=1}^n \overline{d\pi}(x_\ell)}$. Assume $\overline{d\pi}(\sum_{\ell=1}^n x_\ell) = \overline{\sum_{\ell=1}^n \overline{d\pi}(x_\ell)}$ for some $n \in \mathbb{N}$. Let $v \in \mathcal{H}$ and $t \in \mathbb{R}$. By [2, Cor. 3.1.31], the continuity of the map $G \to \mathcal{H}, g \mapsto \pi(g)v$ and since $G$ has the Trotter property we have

$$\begin{aligned}
e^{t\overline{\sum_{\ell=1}^{n+1} \mathrm{d}\pi(x_\ell)}}v &= \lim_{k \to \infty} \left( e^{\frac{t}{k}\overline{\sum_{\ell=1}^{n} \mathrm{d}\pi(x_\ell)}} e^{\frac{t}{k}\overline{\mathrm{d}\pi}(x_{n+1})} \right)^k v \\
&= \lim_{k \to \infty} \left( e^{\frac{t}{k}\overline{\mathrm{d}\pi}(\sum_{\ell=1}^{n} x_\ell)} e^{\frac{t}{k}\overline{\mathrm{d}\pi}(x_{n+1})} \right)^k v \\
&= \pi \left( \lim_{k \to \infty} \left( \exp(\tfrac{t}{k} \sum_{\ell=1}^{n} x_\ell) \exp(\tfrac{t}{k} x_{n+1}) \right)^k \right) v \\
&= \pi \left( \exp(t \sum_{\ell=1}^{n+1} x_\ell) \right) v = e^{t\overline{\mathrm{d}\pi}(\sum_{\ell=1}^{n+1} x_\ell)} v
\end{aligned}$$

Differentiation by $t$ implies $\overline{\sum_{\ell=1}^{n+1} \mathrm{d}\pi(x_\ell)} = \overline{\mathrm{d}\pi}(\sum_{\ell=1}^{n+1} x_\ell)$ so that in particular the domains of these operators coincide. Thus $\mathcal{D}^\infty(N)$ is contained in the domain of each $\overline{\mathrm{d}\pi}(\sum_{\ell=1}^{n} x_\ell)$ and moreover the map

$$\mathrm{d}\pi : \mathfrak{h} \to \mathrm{End}(\mathcal{D}^\infty(N)), \quad x \mapsto \overline{\mathrm{d}\pi}(x)|_{\mathcal{D}^\infty(N)}$$

is linear. Furthermore $\frac{1}{i}\overline{\mathrm{d}\pi}(x)$ is essentially self-adjoint on $\mathcal{D}^\infty(N)$ for every $x \in \mathfrak{h}$ and hence

$$s_\pi(x) = \sup_{\|v\|=1, v \in \mathcal{D}^\infty(N)} \langle i\mathrm{d}\pi(x)v, v \rangle$$

for all $x \in \mathfrak{h}$. Note that for $y = \sum_{\ell=1}^{k} y_\ell \in \mathfrak{h}, y_\ell \in \mathfrak{h}_{j_\ell}$, we have $[y, x_0] \in \mathfrak{h}$ and

$$[\mathrm{d}\pi(y), \mathrm{d}\pi(x_0)] = \sum_{\ell=1}^{k} [\mathrm{d}\pi(y_k), \mathrm{d}\pi(x_0)] = \sum_{\ell=1}^{k} \mathrm{d}\pi([y_k, x_0]) = \mathrm{d}\pi([y, x_0])$$

on $\mathcal{D}^\infty(N)$ since $\mathrm{d}\pi|_{\mathfrak{h}_j}$ is a representation for each $j$. Now we can argue as in the proof of [19, Thm. 3.4], where we take $\mathcal{D}^\infty(N)$ for the space $\mathcal{H}^\infty$ there, to conclude that the map

$$\mathfrak{h} \times \mathcal{D}^\infty(\overline{\mathrm{d}\pi}(x_0)) \to \mathcal{D}^\infty(\overline{\mathrm{d}\pi}(x_0)), \quad (x, v) \mapsto \mathrm{d}\pi(x)v$$

is bilinear and continuous when $\mathcal{D}^\infty(\overline{\mathrm{d}\pi}(x_0))$ is equipped with the $C^\infty$-topology. Thus it extends to a unique continuous bilinear map $\beta : \mathfrak{g} \times \mathcal{D}^\infty(\overline{\mathrm{d}\pi}(x_0)) \to \mathcal{D}^\infty(\overline{\mathrm{d}\pi}(x_0))$. By Lemma 2.7 $\mathcal{D}^\infty(\overline{\mathrm{d}\pi}(x_0)) \subset \mathcal{D}(\mathfrak{g})$ and $\overline{\mathrm{d}\pi}(x)v = \beta(x, v)$. Since $\beta$ takes values in $\mathcal{D}^\infty(\overline{\mathrm{d}\pi}(x_0))$ this implies $\mathcal{D}^\infty(\overline{\mathrm{d}\pi}(x_0)) \subset \mathcal{D}^\infty(\mathfrak{g})$ and as $\beta$ is continuous and bilinear we further conclude that $\mathcal{D}^\infty(\overline{\mathrm{d}\pi}(x_0)) \subset \mathcal{D}_c^\infty(\mathfrak{g})$. This yields (a). If $G$ is locally exponential then it has the Trotter property [18, Prop. 3.5]. Thus (a) and [14, Lem. 3.4] imply that the dense subspace $\mathcal{D}^\infty(\overline{\mathrm{d}\pi}(x_0)) \subset \mathcal{H}$ consists of smooth vectors. Now Lemma 2.6 implies $x_0 \in W_\pi$ and (b) follows. $\square$

**Remark 2.10.** Let $\pi : G \to \mathrm{U}(\mathcal{H})$ be a continuous representation. If $\mathfrak{k} \subset \mathfrak{g}$ is a finite-dimensional subalgebra, then there always exists a corresponding integral subgroup $K \subset G$ ([11, Cor. IV.4.10]) and moreover $\pi|_K$ is smooth since every continuous unitary representation of a finite-dimensional Lie group is smooth. Consider the invariant cone

$$C_\pi := \{x \in \mathfrak{g} : s_\pi(x) \leq 0\}.$$

Then $C_\pi \cap \mathfrak{k} = C_{\pi|_K}$ is a convex $\mathrm{Ad}_K$-invariant cone which is closed in $\mathfrak{k}$ by (1).

## 2.2. Application to oscillator groups.

**Definition 2.11.** (a) Let $(V, \omega)$ be a locally convex symplectic vector space and $\gamma : \mathbb{R} \to \mathrm{Sp}(V, \omega)$ be a one-parameter group of symplectomorphisms defining a smooth action of $\mathbb{R}$ on $V$ and denote by $D := \gamma'(0) : V \to V$ its generator. Then the Lie group

$$G(V, \omega, \gamma) := \mathrm{Heis}(V, \omega) \rtimes_\gamma \mathbb{R}$$

is called an *oscillator group*, where $\mathrm{Heis}(V, \omega) := \mathbb{R} \times_\omega V$ is the *Heisenberg group* with multiplication given by

$$(t, x)(s, y) = \bigl(t + s + \tfrac{1}{2}\omega(x, y), x + y\bigr).$$

The Lie algebra of $G(V, \omega, \gamma)$ is $\mathfrak{g}(V, \omega, \gamma) = \mathrm{heis}(V, \omega) \rtimes_D \mathbb{R}$ with the bracket

$$[(t, v, s), (t', v', s')] = (\omega(v, v'), sDv' - s'Dv, 0).$$

(b) Let $A$ be a self-adjoint operator on a complex Hilbert space $H_A$ with $A \geq 0$ and $\ker A = 0$ and let $\gamma(t) = e^{itA}$ be the corresponding unitary one-parameter group. Equip $V_A := \mathcal{D}^\infty(A)$ with the $C^\infty$-topology generated by the norms $v \mapsto \|A^k v\|, k \in \mathbb{N}_0$. Set $\omega_A(x, y) := \mathrm{Im}\langle Ax, y\rangle$ for $x, y \in V_A$. Then the oscillator group

$$G_A := \mathrm{Heis}(V_A, \omega_A) \rtimes_\gamma \mathbb{R}$$

is called a *standard oscillator group*, cf. [20]. We call $A$ *diagonalizable* if there exists an orthonormal basis $\{e_j : j \in J\}$ of $H_A$ and $\lambda_j \in \mathbb{R}, j \in J$, such that $Ae_j = \lambda_j e_j$ for every $j \in J$.

**Theorem 2.12.** *Let $G_A$ be a standard oscillator group with $A$ diagonalizable and let $\pi : G_A \to \mathrm{U}(\mathcal{H})$ be a continuous positive energy representation, i.e., $\frac{1}{i}\overline{\mathrm{d}\pi}(0, 0, 1) \geq 0$. Then $\pi$ is semibounded and in particular smooth.*

*Proof.* Let $\{e_j : j \in J\}$ be an orthonormal basis of $H_A$ with $Ae_j = \lambda_j e_j, \lambda_j > 0$. Consider the subalgebras $\mathfrak{h}_j := \mathbb{R} \times \mathbb{R}e_j \times \mathbb{R} \subset \mathfrak{g}_A$ and set $\mathfrak{h} := \mathrm{span}_j \mathfrak{h}_j$. Note that each $y = \sum_{k=1}^n y_k \in \mathfrak{h}, y_k \in \mathfrak{h}_{j_k}$, is contained in the finite-dimensional subalgebra $\mathbb{R} \oplus \bigl(\bigoplus_{k=1}^n \mathbb{R}e_{j_k}\bigr) \oplus \mathbb{R} \subset \mathfrak{g}_A$. In view of Remark 2.10 we thus obtain from [24, Prop. 2.17] that

$$s_\pi(t, x, s) = s_\pi\Bigl(t - \frac{\|x\|^2}{2s}, 0, s\Bigr) \quad \text{for all } (t, x, s) \in \mathfrak{h} \text{ with } s \neq 0. \tag{6}$$

By assumption $s_\pi(0, 0, 1) \leq 0$. Choose $j_0 \in J$. Then $C := s_\pi^{-1}(]-\infty, 0]) \cap \mathfrak{h}_{j_0}$ is a convex cone since $\mathfrak{h}_{j_0}$ is a finite-dimensional subalgebra. By (6) the parabola $\bigl\{\bigl(\frac{t^2}{2}, te_{j_0}, 1\bigr) : t \in \mathbb{R}\bigr\}$ is contained in $C$ hence also its convex hull. Since $C$ is a cone we conclude that $s_\pi|_{\mathfrak{h}_{j_0}}$ is bounded in a neighborhood of $x_0 := (1, 0, 1)$. Now (6) implies that $s_\pi|_\mathfrak{h}$ is bounded in a neighborhood of $x_0$. By [18, Prop. 3.14

& Thm. 3.15 & Thm. B.7] $G_A$ has the Trotter property. As $\text{Heis}(V_A)$ is locally exponential, Theorem 2.9(a) and [14, Lem. 3.4] show that $\mathcal{D}^\infty(\overline{\mathbf{d}\pi}(x_0))$ consists of smooth vectors for $\pi|_{\text{Heis}(V_A)}$. Since $\mathcal{D}^\infty(\overline{\mathbf{d}\pi}(x_0)) = \mathcal{D}^\infty(\mathfrak{h}_{j_0})$ by [19, Thm. 3.4] we obtain that

$$G_A \to \mathbb{C}, (t, x, s) \mapsto \langle \pi(t, x, s)v, v \rangle = \langle \pi(0, 0, s)v, \pi(-t, -x)v \rangle$$

is smooth for all $v \in \mathcal{D}^\infty(\overline{\mathbf{d}\pi}(x_0))$. Thus $\pi$ is smooth by [14, Thm. 7.2]. Now Lemma 2.6 implies that $\pi$ is semibounded. $\square$

### 2.3. Application to double extensions of loop groups.

**Definition 2.13.** Let $K$ be a 1-connected simple compact Lie group with Lie algebra $\mathfrak{k}$ and let $\mathcal{L}(K) := C^\infty(S^1, K)$ be the corresponding *loop group*. The Lie algebra of $\mathcal{L}(K)$ is the loop algebra $\mathcal{L}(\mathfrak{k}) := C^\infty(S^1, \mathfrak{k})$. Since $K$ is compact we may choose an invariant symmetric positive definite form $\langle \cdot, \cdot \rangle$ on $\mathfrak{k}$ which is normalized in the sense of [17, Def. 3.3]. This yields the double extension

$$\widehat{\mathcal{L}}(\mathfrak{k}) := (\mathbb{R} \oplus_\omega \mathcal{L}(\mathfrak{k})) \rtimes_D \mathbb{R},$$

where $\omega(\xi, \eta) = \frac{1}{2\pi} \int_0^{2\pi} \langle \xi'(t), \eta(t) \rangle dt$ and $D\xi = \xi'$, see [17, Exa. 2.4]. Note that $\mathfrak{t} := \mathbb{R} \oplus \{0\} \oplus \mathbb{R} \subset \widehat{\mathcal{L}}(\mathfrak{k})$ is an abelian subalgebra. The Lie algebra $\mathbb{R} \oplus_\omega \mathcal{L}(\mathfrak{k})$ integrates to a 2-connected Lie group $\widetilde{\mathcal{L}}(K)$ ([17, Thm. 3.4]). The rotation action $\alpha$ of $\mathbb{R}$ on $\mathcal{L}(K)$ lifts uniquely to a smooth action on $\widetilde{\mathcal{L}}(K)$ which yields a 2-connected Fréchet-Lie group

$$\widehat{\mathcal{L}}(K) := \widetilde{\mathcal{L}}(K) \rtimes_\alpha \mathbb{R}$$

with Lie algebra $\widehat{\mathcal{L}}(\mathfrak{k})$ (cf. [17, Def. 3.5]).

A continuous representation $\pi : \widehat{\mathcal{L}}(K) \to \mathrm{U}(\mathcal{H})$ is called a *positive energy representation* if the self-adjoint generator $-i\overline{\mathbf{d}\pi}(0, 0, 1)$ is non-negative.

**Remark 2.14.** (a) Using the identification $S^1 = \mathbb{R}/2\pi\mathbb{Z}$ we may identify

$$\mathcal{L}(K) = \{f \in C^\infty(\mathbb{R}, K) : (\forall t \in \mathbb{R})\ f(t + 2\pi) = f(t)\},$$

$$\mathcal{L}(\mathfrak{k}) = \{\xi \in C^\infty(\mathbb{R}, \mathfrak{k}) : (\forall t \in \mathbb{R})\ \xi(t + 2\pi) = \xi(t)\}.$$

(b) The adjoint action of $\mathcal{L}(K)$ on the double extension $\mathfrak{g} = \widehat{\mathcal{L}}(\mathfrak{k})$ is given by

$$\text{Ad}_\mathfrak{g}(g)(z, \xi, t) = \left(z + \langle \delta^l(g), \xi \rangle - \frac{t}{2}\|\delta^r(g)\|^2, \text{Ad}(g)\xi - t\delta^r(g), t\right) \quad (7)$$

where $\delta^r(g) = g'g^{-1} \in \mathcal{L}(\mathfrak{k})$ denotes the *right logarithmic derivative* and $\delta^l(g) = g^{-1}g' \in \mathcal{L}(\mathfrak{k})$ denotes the *left logarithmic derivative* of $g \in \mathcal{L}(K)$, cf. [17, Sect. 2.3].

**Proposition 2.15.** *Let $\pi : \widehat{\mathcal{L}}(K) \to \mathrm{U}(\mathcal{H})$ be a continuous positive energy representation. Then $W_\pi \cap \mathfrak{t} \neq \emptyset$. In particular, $W_\pi$ is non-empty.*

*Proof.* Set $C_\pi := \{x \in \widehat{\mathcal{L}}(\mathfrak{k}) : s_\pi(x) \leq 0\}$ and $\mathfrak{g} := \widehat{\mathcal{L}}(\mathfrak{k})$. Consider the subalgebra $\widehat{\mathfrak{k}} := \mathbb{R} \oplus \mathfrak{k} \oplus \mathbb{R}$ of $\widehat{\mathcal{L}}(\mathfrak{k})$ where we identify $\mathfrak{k}$ with the constant loops in $\mathfrak{k}$. Since $K$ is compact there exists a basis $(e_i)_{i=1}^n$ of $\mathfrak{k}$ such that $\exp(2\pi e_i) = \mathbf{1}$ for all $i \in \{1, \ldots, n\}$. Indeed, first choose a basis in a subalgebra $\mathfrak{t}_\mathfrak{k} \subset \mathfrak{k}$, where $\mathfrak{t}_\mathfrak{k}$ is the Lie algebra of a maximal torus, and then use $\mathrm{Ad}(K)\mathfrak{t}_\mathfrak{k} = \mathfrak{k}$. Since $\widehat{\mathfrak{k}}$ is finite-dimensional, $C := C_\pi \cap \widehat{\mathfrak{k}}$ is a closed convex invariant cone in $\widehat{\mathfrak{k}}$ which contains $(0, 0, 1)$ by assumption. By the choice of $e_i$ we may consider $\nu_{i,j}(t) := \exp(tje_i)$ as elements $\nu_{i,j} \in C^\infty(S^1, K)$ for all $i \in \{1, \ldots, n\}, j \in \mathbb{Z}$. As $C_\pi$ is $\mathrm{Ad}_\mathfrak{g}$-invariant, equation (7) yields

$$\mathrm{Ad}(\nu_{i,j})(0, 0, 1) = \left(-\tfrac{1}{2}\|je_i\|^2, -je_i, 1\right) \in C. \tag{8}$$

It is easy to verify that [2]

$$\left\{\left(-\tfrac{1}{2}(\lambda^2 + \tfrac{1}{4})\|e_i\|^2, \lambda e_i, 1\right) : \lambda \in \mathbb{R}\right\} \subset \mathrm{conv}\left\{\left(-\tfrac{1}{2}j^2\|e_i\|^2, je_i, 1\right) : j \in \mathbb{Z}\right\}$$

for all $i = 1, \ldots, n$. Thus equation (8) and the convexity of $C$ imply

$$\left\{\left(-\tfrac{n}{2}\|x\|_b^2 - d, x, 1\right) : x \in \mathfrak{k}, d \geq c_0\right\} \subset C,$$

where $\|\sum_{i=1}^n \lambda_i e_i\|_b := (\sum_{i=1}^n \lambda_i^2 \|e_i\|^2)^{1/2}$ and $c_0 := \tfrac{1}{8n}\sum_{i=1}^n \|e_i\|^2$. Hence there exists $c_1 > 0$ such that $\{(-c_1\|x\|^2 - d, x, 1) : x \in \mathfrak{k}, d \geq c_0\} \subset C$. Since $C$ is a cone this entails an element $x_0 = (z_0, 0, 1) \in \widehat{\mathfrak{k}}$ which is contained in the interior $C^0$ of $C$ in $\widehat{\mathfrak{k}}$.

Now choose an open 0-neighborhood $U_\mathfrak{k} \subset \mathfrak{k}$ such that $\exp : U_\mathfrak{k} \to U_K$ is a diffeomorphism onto an open $\mathbf{1}$-neighborhood $U_K \subset K$. Since $K$ is compact, for every $\xi \in C^\infty(\mathbb{R}, \mathfrak{k})$ there exists a unique $\gamma_\xi \in C^\infty(\mathbb{R}, K)$ such that $\gamma_\xi(0) = \mathbf{1}$ and $\delta^r(\gamma_\xi) = \xi$. Moreover there exists an open convex symmetric 0-neighborhood $\widetilde{U} \subset \mathfrak{k}$ such that the map

$$\mathrm{Evol} : C^\infty(S^1, \widetilde{U}) \to C^\infty([0, 2\pi], U_K), \quad \xi \mapsto \gamma_\xi|_{[0, 2\pi]}$$

is well-defined. Note that Evol is smooth due to the smooth dependence of solutions of ordinary differential equations on (time-dependent) vector fields. Set $U_{\mathcal{L}(\mathfrak{k})} := C^\infty(S^1, \tfrac{1}{2}\widetilde{U})$ which is an open subset of $\mathcal{L}(\mathfrak{k})$. For $\xi \in U_{\mathcal{L}(\mathfrak{k})}$ and $t > \tfrac{1}{2}$ define $g_{\xi,t} := \gamma_{\tfrac{-1}{t}\xi}$, $y_\xi := \tfrac{1}{2\pi}\exp^{-1}(g_{\xi,t}(2\pi))$ and $h_\xi(s) := \exp(sy_\xi)$. Then we may consider $h_\xi g_{\xi,t}^{-1}$ as element of $\mathcal{L}(K)$. The map

$$F : \mathbb{R} \times U_{\mathcal{L}(\mathfrak{k})} \times \,]\tfrac{1}{2}, \infty[ \,\to \widehat{\mathcal{L}}(\mathfrak{k}), \quad F(z, \xi, t) = \mathrm{Ad}_\mathfrak{g}(h_\xi g_{\xi,t}^{-1})(z, \xi, t),$$

is continuous (even smooth) since Evol is smooth. For $\xi \in U_{\mathcal{L}(\mathfrak{k})}, t > 0$ we calculate with [11, Rem. II.3.4]

$$\mathrm{Ad}(h_\xi g_{\xi,t}^{-1})\xi - t\delta^r(h_\xi g_{\xi,t}^{-1}) = \mathrm{Ad}(h_\xi g_{\xi,t}^{-1})\xi - t(\delta^r(h_\xi) + \mathrm{Ad}(h_\xi)\delta^r(g_{\xi,t}^{-1}))$$
$$= \mathrm{Ad}(h_\xi)\mathrm{Ad}(g_{\xi,t}^{-1})\xi - ty_\xi + t\mathrm{Ad}(h_\xi)\mathrm{Ad}(g_{\xi,t}^{-1})\delta^r(g_{\xi,t})$$
$$= -ty_\xi$$

---

[2] This follows from the geometry of the standard parabola: If $P_d(t) := (\tfrac{1}{2}t^2 + d, t)$ and $c > 0$ then $P_{\tfrac{1}{8}c^2}(\mathbb{R}) \subset \mathrm{conv}(P_0(\mathbb{Z}c))$, where $\mathrm{conv}(M)$ denotes the convex hull of $M \subset \mathbb{R}^2$.

since $\delta^r(g_{\xi,t}) = \frac{-1}{t}\xi$. Now (7) entails that $F$ takes values in $\widehat{\mathfrak{k}}$. Recall the element $x_0 = (z_0, 0, 1) \in C^0$ from above. Since $F(x_0) = x_0$ and $F$ is continuous there exists an open $x_0$-neighborhood $V$ in $\widehat{\mathcal{L}}(\mathfrak{k})$ such that $F(V) \subset C = C_\pi \cap \widehat{\mathfrak{k}}$. As $C_\pi$ is $\mathrm{Ad}_{\mathfrak{g}}$-invariant we conclude that $V \subset C_\pi$. Hence $x_0 \in W_\pi \cap \mathfrak{t}$. □

**Theorem 2.16.** *Let $K$ be a 1-connected simple compact Lie group and $\widehat{\mathcal{L}}(K)$ be the corresponding double extension of the loop group $\mathcal{L}(K)$. Then every continuous positive energy representation $\pi: \widehat{\mathcal{L}}(K) \to \mathrm{U}(\mathcal{H})$ is semibounded and in particular smooth.*

*Proof.* Let $G := \widehat{\mathcal{L}}(K), \mathfrak{g} := \widehat{\mathcal{L}}(\mathfrak{k})$, choose a basis $(b_j)_{1 \leq j \leq N}$ of $\mathfrak{k}$ and set

$$e_{n,j} \in \mathcal{L}(\mathfrak{k}), e_{n,j}(t) := \sin(nt)b_j, \quad f_{n,j} \in \mathcal{L}(\mathfrak{k}), f_{n,j}(t) := \cos(nt)b_j$$

for $n \in \mathbb{N}_0, 1 \leq j \leq N$. Then each $\mathfrak{g}_{n,j} := \mathbb{R} \oplus_\omega (\mathbb{R}e_{n,j} + \mathbb{R}f_{n,j}) \rtimes_D \mathbb{R}$ is a subalgebra of $\mathfrak{g}$ and $\mathrm{span}_{n,j}\mathfrak{g}_{n,j}$ is dense in $\mathfrak{g}$. By Proposition 2.15 there exists an element $x_0 \in W_\pi \cap \mathfrak{t} \subset \bigcap_{n,j} \mathfrak{g}_{n,j}$. Since $\widehat{\mathcal{L}}(K)$ satisfies the Trotter property ([6]) Theorem 2.9 yields $\mathcal{D}^\infty(\overline{\mathrm{d}\pi}(x_0)) = \mathcal{D}_c^\infty(\mathfrak{g})$. Since $\widetilde{\mathcal{L}}(K)$ is locally exponential, $\mathcal{D}^\infty(\overline{\mathrm{d}\pi}(x_0))$ consists of smooth vectors for $\pi|_{\widetilde{\mathcal{L}}(K)}$ by [14, Lem. 3.4]. Let $T \subset G$ be the integral subgroup corresponding to $\mathfrak{t} \subset \mathfrak{g}$. Then $\mathcal{D}^\infty(\overline{\mathrm{d}\pi}(x_0)) = \mathcal{H}^\infty(\pi|_T)$ by [19, Thm. 3.4]. Therefore

$$G = \widetilde{\mathcal{L}}(K) \rtimes \mathbb{R} \to \mathbb{C}, \quad (x,t) \mapsto \langle \pi(x,t)v, v \rangle = \langle \pi(t)v, \pi(x)^{-1}v \rangle$$

is smooth for every $v \in \mathcal{D}^\infty(\overline{\mathrm{d}\pi}(x_0))$. Now [14, Thm. 7.2] entails that $\pi$ is smooth with $\mathcal{H}^\infty(\pi) = \mathcal{D}^\infty(\overline{\mathrm{d}\pi}(x_0))$. Thus $\pi$ is semibounded as $W_\pi \neq \emptyset$. □

## 2.4. Application to the Virasoro group.

**Definition 2.17.** The Lie algebra $\mathcal{V}(S^1)$ of smooth vector fields on the unit circle $S^1$ is the Lie algebra of the group of orientation preserving diffeomorphisms $\mathrm{Diff}(S^1)_+$ of the unit circle. Let $\partial_\theta := \frac{d}{d\theta}$ denote the vector field generating the rigid rotations of $S^1$. Note that $\mathcal{V}(S^1) = C^\infty(S^1)\partial_\theta$ and

$$[f\partial_\theta, g\partial_\theta] = (fg' - gf')\partial_\theta, \quad f, g \in C^\infty(S^1).$$

The *Virasoro algebra* is the central extension $\mathfrak{vir} := \mathbb{R} \oplus_\omega \mathcal{V}(S^1)$ of $\mathcal{V}(S^1)$ by the cocycle

$$\omega(f\partial_\theta, g\partial_\theta) := \int_0^{2\pi}(f''' + f')g d\theta,$$

cf. [15, Sect. 8.2]. The *Virasoro group* Vir is the (up to isomorphism unique) simply connected regular Lie group with Lie algebra $\mathfrak{vir}$. Vir is a central extension of the universal covering group of $\mathrm{Diff}(S^1)_+$. A continuous unitary representation $\pi: \mathrm{Vir} \to \mathrm{U}(\mathcal{H})$ is called a *positive energy representation* if $\frac{1}{i}\overline{\mathrm{d}\pi}(0, \partial_\theta) \geq 0$.

**Theorem 2.18.** *Every continuous positive energy representation $\pi: \mathrm{Vir} \to \mathrm{U}(\mathcal{H})$ of the Virasoro group is smooth and semibounded with $(1, \partial_\theta) \in W_\pi$.*

*Proof.* Set $G := \mathrm{Vir}, \mathfrak{g} := \mathfrak{vir}$ and $C_\pi := \{x \in \mathfrak{g} : s_\pi(x) \leq 0\}$. For $n \in \mathbb{N}$ set

$$e_n \in \mathcal{V}(S^1), e_n(\theta) := \sin(n\theta)\partial_\theta, \quad f_n \in \mathcal{V}(S^1), f_n(\theta) := \cos(n\theta)\partial_\theta.$$

Then each $\mathfrak{g}_n := \mathbb{R} \oplus_\omega (\mathbb{R}e_n + \mathbb{R}f_n + \mathbb{R}\partial_\theta)$ is a subalgebra of $\mathfrak{g}$ and $\mathrm{span}_n \mathfrak{g}_n$ is dense in $\mathfrak{g}$. Let $\mathfrak{t} := \mathbb{R} \oplus \mathbb{R}\partial_\theta$ and denote by $T_r := \exp(\mathbb{R}\partial_\theta)$ the subgroup of rigid rotations. Consider the projection

$$p_\mathfrak{t} : \mathfrak{g} \to \mathfrak{t}, \quad (z, f\partial_\theta) \mapsto \left(z, \int_{T_r} \mathrm{Ad}(\varphi)(f\partial_\theta) d\mu_{T_r}(\varphi)\right) = \left(z, \frac{1}{2\pi}\int_0^{2\pi} f(s)ds \cdot \partial_\theta\right)$$

where $\mu_{T_r}$ denotes the normalized Haar measure on $T_r$. Set $d := (0, \partial_\theta)$ and denote by $\mathcal{O}_d^{\mathfrak{g}_n}$ the adjoint orbit of $d$ in $\mathfrak{g}_n$. From the proof of [15, Thm. 8.14] we obtain

$$d + \left(1, \tfrac{1}{\pi(n^2-1)}\partial_\theta\right) \in p_\mathfrak{t}(\mathcal{O}_d^{\mathfrak{g}_n}) \subset \overline{\mathrm{conv}}(\mathcal{O}_d^{\mathfrak{g}_n})$$

for $n \geq 2$, where $\overline{\mathrm{conv}}(\mathcal{O}_d^{\mathfrak{g}_n})$ denotes the closed convex hull of $\mathcal{O}_d^{\mathfrak{g}_n}$ in $\mathfrak{g}_n$. Since $d \in C_\pi$ and $C_\pi \cap \mathfrak{g}_n$ is a closed convex invariant cone in $\mathfrak{g}_n$ (Remark 2.10) we obtain $d + (1, \frac{1}{\pi(n^2-1)}\partial_\theta) \in C_\pi \cap \mathfrak{t}$ for $n \geq 2$. As $d \in C_\pi$ and $C_\pi \cap \mathfrak{t}$ is a closed convex cone in $\mathfrak{t}$ we conclude

$$\mathbb{R}_+(1,0) + \mathbb{R}_+(0, \partial_\theta) \subset C_\pi \cap \mathfrak{t}, \quad \text{for } \mathbb{R}_+ := \{t \in \mathbb{R} : t \geq 0\}.$$

In particular, $x_0 := (1, \partial_\theta)$ lies in the interior of $C_\pi \cap \mathfrak{t}$ in $\mathfrak{t}$. Let $U$ be an open neighborhood of $x_0$ in $\mathfrak{t}$ with $U \subset C_\pi \cap \mathfrak{t}$. Set $W_{\max} := \{(z, f) \in \mathfrak{g} : f > 0\}$. According to [15, Prop. 8.12] there is a continuous map

$$F : W_{\max} \to \mathfrak{t}, \quad (z, f) \mapsto (\beta(z, f), \alpha(z, f))$$

with $F(x_0) = x_0$ such that $F(z, f)$ and $(z, f)$ lie on the same $\mathrm{Ad}_G$-orbit for every $(z, f) \in W_{\max}$. As $C_\pi$ is $\mathrm{Ad}_G$-invariant we conclude that $x_0 \in F^{-1}(U) \subset C_\pi$. By the continuity of $F$ the set $F^{-1}(U)$ is an open $x_0$-neighborhood in $\mathfrak{g}$. Therefore $x_0 \in W_\pi$. Since Vir has the Trotter property ([18, Cor. 3.17]), Theorem 2.9 yields $\mathcal{D}^\infty(\overline{d\pi}(x_0)) = \mathcal{D}_c^\infty(\mathfrak{g})$. The Virasoro group is $R$-regular by [6]. Thus Proposition 2.3 entails $\mathcal{D}^\infty(\overline{d\pi}(x_0)) = \mathcal{D}_c^\infty(\mathfrak{g}) = \mathcal{H}^\infty$. In particular, $\pi$ is smooth and thus semibounded since $x_0 \in W_\pi$. □

**Remark 2.19.** Let $\pi : \mathrm{Vir} \to \mathrm{U}(\mathcal{H})$ be a continuous positive energy representation and $d = (0, \partial_\theta), \mathfrak{t} = \mathbb{R} \oplus \mathbb{R}\partial_\theta$. In [7, Chapter 1] it is shown that $\mathcal{D}_c^\infty(\mathfrak{vir}) = \mathcal{D}^\infty(\mathfrak{t})$ if $\frac{1}{i}\overline{d\pi}(1, 0) = c\mathbf{1}, c \in \mathbb{R}$, and $\frac{1}{i}\overline{d\pi}(d)$ is diagonalizable with finite-dimensional eigenspaces, see also [5, Sect. 4], [3, App. A] and [4, Thm. 3.4] where the eigenspaces of $\frac{1}{i}\overline{d\pi}(d)$ are not assumed to be finite-dimensional.

## 3. Analytic vectors

In this section $G$ denotes a Banach–Lie group. Then we may choose an open 0-neighborhood $\widetilde{U} \subset \mathfrak{g}$ such that $\exp|_{\widetilde{U}} : \widetilde{U} \to \exp(\widetilde{U})$ is a diffeomorphism and the map

$$m : \widetilde{U} \times \widetilde{U} \to \mathfrak{g}, \quad (x, y) \to \exp^{-1}(\exp(x)\exp(y))$$

is analytic, cf. [11, Cor. IV.1.10 & Def. IV.1.5]. Hence there exists an open 0-neighborhood $U \subset \widetilde{U}$ and a holomorphic map $m_{\mathbb{C}} : U_{\mathbb{C}} \times U_{\mathbb{C}} \to \mathfrak{g}_{\mathbb{C}}$, where $U_{\mathbb{C}} := U + iU$, such that $m_{\mathbb{C}}$ extends $m|_U$ and is given by the Hausdorff series. We write $x * y := m_{\mathbb{C}}(x, y)$ for $x, y \in U_{\mathbb{C}}$. Consider the map

$$\psi : U \times U \to \mathfrak{g}_{\mathbb{C}}, \quad (x, y) \to x * iy.$$

Since $\mathtt{d}\psi(0,0)(x,y) = x+iy$ we may assume by the Inverse Function Theorem (after shrinking $U$) that $\psi$ is a diffeomorphism onto an open 0-neighborhood $\psi(U \times U) \subset \mathfrak{g}_{\mathbb{C}}$. By shrinking $U$ further to a convex 0-neighborhood we may assume that three-fold products are defined and $(a * b) * c = a * (b * c)$ holds for all $a, b, c \in U_{\mathbb{C}}$.

**Remark 3.1.** (a) Since $(b*a)*(-a) = b$ and $(-a)*(a*b) = b$ for all $a, b, c \in U$ these equations hold by analytic continuation also for all $a, b, c \in U_{\mathbb{C}}$. With $\ell_a(b) := a*b$ and $r_a(b) := b*a$ we conclude that $\mathtt{d}\ell_a(b)$ and $\mathtt{d}r_a(b)$ are $\mathbb{C}$-linear isomorphisms for all $a, b \in U_{\mathbb{C}}$.

(b) Let $I \subset \mathbb{R}$ be an interval and $\gamma : I \to U_{\mathbb{C}}$ be a $C^1$ curve. Then the *(right) logarithmic derivative* of $\gamma$ is defined by

$$\delta^r(\gamma)_t := \mathtt{d}r_{-\gamma(t)}(\gamma(t))(\gamma'(t)) = (\mathtt{d}r_{\gamma(t)}(0))^{-1}(\gamma'(t)),$$

see also [11, Sect. II.4]. By [9, Lemma 13] $\delta^r(\gamma)_t = \int_0^1 e^{\mathrm{sad}\gamma(t)}\gamma'(t)ds$.

Now let $\pi : G \to \mathrm{U}(\mathcal{H})$ be a semibounded unitary representation.

**Definition 3.2.** Define

$$\widetilde{\rho} : \mathfrak{g} \times W_\pi \to B(\mathcal{H}), \quad (x, w) \mapsto \pi(\exp(x))e^{i\mathtt{d}\pi(w)}$$

and set $V_\pi := \psi(U \times (U \cap W_\pi))$ and $\rho : V_\pi \to B(\mathcal{H}), \rho(z) = \widetilde{\rho}(\psi^{-1}(z))$.

We want to show that $\rho$ is holomorphic as a map from the open subset $V_\pi \subset \mathfrak{g}_{\mathbb{C}}$ to the bounded operators $B(\mathcal{H})$ (equipped with the operator norm topology).

**Remark 3.3.** For $v \in \mathcal{H}^\infty$ the map $\mathfrak{g} \times W_\pi \to \mathcal{H}, (x, w) \to \widetilde{\rho}(x,w)v$ is $C^1$ by [10, Lem. 5.2]. Thus $V_\pi \to \mathcal{H}, z \mapsto \rho(z)v$ is $C^1$ for all $v \in \mathcal{H}^\infty$. Moreover $e^{i\mathtt{d}\pi(w)}\mathcal{H}^\infty \subset \mathcal{H}^\infty$ for all $w \in W_\pi$, by the proof of [10, Lem. 5.2].

**Lemma 3.4.** Let $y \in U \cap W_\pi$ and $U_y \subset U$ be a convex 0-neighborhood such that $(iU_y) * iy \subset V_\pi = U * (i(U \cap W_\pi))$. Let $x \in U_y$. For every $t \in [0,1]$ there exist unique elements $a(t) \in U, w(t) \in U \cap W_\pi$ such that $(tix)*iy = a(t)*iw(t)$. Then

$$e^{ti\mathtt{d}\pi(x)}e^{i\mathtt{d}\pi(y)} = \pi(\exp(a(t)))e^{i\mathtt{d}\pi(w(t))} \quad \forall t \in [0,1]. \tag{9}$$

*Proof.* Note that for $y \in U \cap W_\pi$ it is always possible to choose $U_y$ with the stated properties since $U \to \mathfrak{g}_{\mathbb{C}}, z \mapsto iz*iy$ is continuous. The first statement follows from the fact that $\psi : U \times U \to U*iU$ is a diffeomorphism. Thus it remains to prove (9). Here we argue in a similar fashion to the proof of [10, Thm. 5.4]. Consider the map $\gamma : [0,1] \to \mathcal{H}, \gamma(t) := \pi(\exp(a(t)))e^{i\mathtt{d}\pi(w(t))}v$ which is $C^1$ by Remark 3.3. Note $\frac{d}{ds}\big|_{s=t}\big(a(s)*(-a(t))\big) = \mathrm{d}r_{-a(t)}(a(t))(a'(t)) = \delta^r(a)_t$ and denote the complex-linear extension of $\mathrm{d}\pi$ again by $\mathrm{d}\pi : \mathfrak{g}_{\mathbb{C}} \to \mathrm{End}(\mathcal{H}^\infty)$. Then [10, Lem. 5.2 & Lem. 5.5] yields

$$\begin{aligned}\gamma'(t) &= \mathrm{d}\pi\big(\delta^r(a)_t\big)\gamma(t) + \pi\big(\exp(a(t))\big)\mathrm{d}\pi\Big(\int_0^1 e^{is\cdot\mathrm{ad}w(t)}iw'(t)ds\Big)e^{i\mathrm{d}\pi(w(t))}v \\ &= \mathrm{d}\pi(\delta^r(a)_t)\gamma(t) + \pi(\exp(a(t)))\mathrm{d}\pi(\delta^r(iw)_t)e^{i\mathrm{d}\pi(w(t))}v \\ &= \mathrm{d}\pi\big(\delta^r(a)_t + e^{\mathrm{ad}a(t)}\delta^r(iw)_t\big)\gamma(t) \\ &= \mathrm{d}\pi\big(\delta^r(a*w)_t\big)\gamma(t) \qquad \text{(by [9, Lem. 12])} \\ &= i\mathrm{d}\pi(x)\gamma(t),\end{aligned}$$

where the last equation follows from $(tix)*iy = a(t)*iw(t)$ and [9, Lem. 11]. Since $\gamma(0) = e^{i\mathrm{d}\pi(y)}$ the assertion now follows from [8, Thm. 1.5]. $\square$

**Proposition 3.5.** *The map* $\rho : V_\pi \to B(\mathcal{H})$ *is holomorphic.*

*Proof.* Let $v \in \mathcal{H}^\infty$ and consider the map $\rho^v : V_\pi \to \mathcal{H}, z \mapsto \rho(z)v$ which is $C^1$ by Remark 3.3. Let $y \in U \cap W_\pi$ and $x \in \mathfrak{g}$. Then for small $t \in \mathbb{R}$ we have $\rho((tix)*(iy)) = e^{ti\mathrm{d}\pi(x)}e^{i\mathrm{d}\pi(y)}$ by Lemma 3.4 and $\rho((tx)*(iy)) = e^{t\mathrm{d}\pi(x)}e^{i\mathrm{d}\pi(y)}$ by the definition of $\rho$. This implies that $\mathrm{d}\rho^v(iy)$ is $\mathbb{C}$-linear since $\mathrm{d}r_{iy}(0)$ is a $\mathbb{C}$-linear isomorphism. For $x \in U, y \in U \cap W_\pi$ and $z \in \mathfrak{g}_{\mathbb{C}}$

$$\rho(x*(tz)*iy) = \pi(\exp(x))\rho((tz)*iy)$$

holds for small $t \in \mathbb{R}$. Thus Remark 3.1(a) implies that $\rho^v$ is complex-differentiable at every point $x*iy \in V_\pi$. Since $\rho^v$ is $C^1$ we conclude that it is holomorphic. As $\pi$ is semibounded and $\|\widetilde{\rho}(x,w)\| = \|e^{i\mathrm{d}\pi(w)}\| = e^{s_\pi(w)}$, the map $\rho = \widetilde{\rho} \circ \psi^{-1}$ is locally bounded. Hence $\rho$ is holomorphic by [12, Lem. 3.4] since $\mathcal{H}^\infty \subset \mathcal{H}$ is dense. $\square$

Recall that Banach–Lie groups are analytic.

**Definition 3.6.** Let $\pi : G \to \mathrm{U}(\mathcal{H})$ be a continuous representation of the Banach–Lie group $G$. A vector $v \in \mathcal{H}$ is called *analytic* if the orbit map $G \to \mathcal{H}, g \mapsto \pi(g)v$ is analytic, see [16].

**Theorem 3.7.** *Let $G$ be a Banach–Lie group and $\pi : G \to \mathrm{U}(\mathcal{H})$ be a semibounded representation. Then for every $x_0 \in W_\pi$ the map $G \to B(\mathcal{H}), g \mapsto \pi(g)e^{i\mathrm{d}\pi(x_0)}$ is analytic. In particular, there exists a dense subspace of analytic vectors for $\pi$.*

*Proof.* Let $U \subset \mathfrak{g}$ be as above. Since $e^{i\mathrm{d}\pi(x_0)} = (e^{i\mathrm{d}\pi(x_0/n)})^n$ we may assume w.l.o.g. that $x_0 \in W_\pi \cap U$. Then the map

$$U \to \mathcal{H}, \quad x \mapsto \pi(\exp(x))e^{i\mathrm{d}\pi(x_0)} = \rho(x*ix_0)$$

is analytic since $x \to x*ix_0$ is analytic and $\rho$ is holomorphic (Proposition 3.5). Hence $G \to B(\mathcal{H}), g \mapsto \pi(g)e^{i\mathrm{d}\pi(x_0)}$ is analytic in a **1**-neighborhood and therefore analytic since $\pi$ is a representation. As $e^{i\mathrm{d}\pi(x_0)}\mathcal{H} \subset \mathcal{H}$ is dense ($e^{i\mathrm{d}\pi(x_0)}$ is self-adjoint and injective), we conclude that $e^{i\mathrm{d}\pi(x_0)}\mathcal{H}$ is a dense subspace consisting of analytic vectors. $\square$

**Remark 3.8.** Let $\pi : G \to \mathrm{U}(\mathcal{H})$ be semibounded and $x_0 \in W_\pi$. Theorem 3.7 shows in particular that the map $G \to B(\mathcal{H}), g \mapsto \pi(g)e^{i\mathtt{d}\pi(x_0)}$ is smooth if $G$ is a Banach–Lie group. Using different techniques this was also obtained in [19] in the more general setting when $G$ is only assumed to be metrizable.

## 4. Open problems

If $\pi : G \to \mathrm{U}(\mathcal{H})$ is a continuous representation with $W_\pi \neq \emptyset$ then Theorem 2.9 shows the existence of a dense space of smooth vectors for $\pi$ under certain conditions, including conditions on the structure of the Lie algebra $\mathfrak{g}$. Though this theorem applies to many important examples of infinite-dimensional Lie groups, a more general result would be desirable.

**Conjecture 4.1.** *Let $G$ be a Lie group with an exponential map and $\pi : G \to \mathrm{U}(\mathcal{H})$ a continuous representation with $W_\pi \neq \emptyset$. Then $\pi$ is semibounded, resp., smooth.*

Here the crucial part is to show that $\pi$ is smooth. If $\pi$ is semibounded then $\mathcal{H}^\infty(\pi) = \mathcal{D}^\infty(\overline{\mathtt{d}\pi}(x_0))$ for any $x_0 \in W_\pi$ ([19]). Thus in the situation of the preceding conjecture one can pick $x_0 \in W_\pi$ and take the dense subspace $\mathcal{D}^\infty(\overline{\mathtt{d}\pi}(x_0))$ as a natural candidate for $\mathcal{H}^\infty(\pi)$. Nevertheless, a proof of the conjecture seems to require new techniques as for $x, y \in \mathfrak{g}$ it is not even clear that $\overline{\mathtt{d}\pi}(x)$ and $\overline{\mathtt{d}\pi}(y)$ are defined on a common dense domain.

According to Theorem 3.7 every semibounded representation $\pi$ of a Banach–Lie group $G$ has a dense space of analytic vectors, where analytic vectors are understood in the sense of Definition 3.6. Even though locally convex Lie groups are in general not analytic, there are natural notions of analytic vectors for their representations, cf. [16]. A generalization of Theorem 3.7 from the Banach case to general locally convex Lie groups would be desirable. As $\mathcal{H}^\infty = \mathcal{D}^\infty(\overline{\mathtt{d}\pi}(x_0)), x_0 \in W_\pi$, for a semibounded representation $\pi$ it seems natural that the space of analytic vectors $\mathcal{H}^\omega$ of $\pi$ is also determined by $\mathtt{d}\pi(x_0)$ for $x_0 \in W_\pi$:

**Problem 4.2.** *Let $G$ be a Lie group with an exponential map, $\pi : G \to \mathrm{U}(\mathcal{H})$ a semibounded representation and $x_0 \in W_\pi$. Show that $\mathcal{H}^\omega = \mathcal{H}^\omega(\mathtt{d}\pi(x_0))$, where*

$$\mathcal{H}^\omega = \Big\{v \in \mathcal{H}^\infty : \sum_{n=0}^\infty \tfrac{\|\mathtt{d}\pi(x)^n v\|}{n!} < \infty \text{ for } x \text{ in a neighbh. of } 0\Big\},$$
$$\mathcal{H}^\omega(\mathtt{d}\pi(x_0)) = \Big\{v \in \mathcal{H}^\infty : (\exists t > 0) \sum_{n=0}^\infty \tfrac{t^n \|\mathtt{d}\pi(x_0)^n v\|}{n!} < \infty\Big\}.$$

Christoph Zellner, Department of Mathematics, Friedrich-Alexander-University, Erlangen-Nuremberg, Cauerstrasse 11, 91058 Erlangen, Germany
E-mail: zellner@mi.uni-erlangen.de